\documentclass[12pt]{article}
\usepackage{amssymb}
\usepackage{color}
\usepackage{indentfirst}

\newtheorem{theorem}{Theorem}
\newtheorem{lemma}{Lemma}

\newtheorem{example}{Example}

\def\qi#1 {\fbox {\footnote {\ }}\ \footnotetext { From Qi: {\color{red}#1}}}

\begin{document}
\title{Infinite families of $2$-designs from a class of cyclic codes with two non-zeros}
\author{Xiaoni Du$^*$,  Rong Wang$^*$£¬ Cuiling Fan$^{*\dag}$ \\
$^{*}$College of Mathematics and Statistics, Northwest Normal \\University,
 Lanzhou,  Gansu 730070,  P.R. China\\
Email: \{ymldxn, rongw113\}@126.com\\
$^{*\dag}$ School of Mathematics,
 Southwest Jiaotong \\
 University, Chengdu, China\\
 Email: cuilingfan@163.com}
\date{}
\maketitle

\begin{abstract}
Combinatorial $t$-designs have wide applications in coding theory, cryptography, communications and statistics. It is well known that the supports of
all codewords with a fixed weight in a
code may give  a $t$-design. In this paper, we first determine the weight distribution of a class of linear codes derived from the dual of extended
cyclic code with two non-zeros. We then obtain infinite families of $2$-designs  and explicitly compute their parameters from the supports of all the
codewords with a fixed weight in the codes. By simple counting argument, we obtain exponentially many $2$-designs.

\end{abstract}

\section{Introduction}
Let $k,t$ and $v$ are positive integers with $1 \leq t \leq k
\leq v.$ Let $\mathcal{P}$ be a set of $v$ elements and $\mathcal{B}$ be a set of $k$-subsets of $\mathcal{P}$.  If every $t$-subset of $\mathcal{P}$ is
contained in exactly $\lambda$ elements of $\mathcal{B},$
then we call the pair $\mathbb{D}=(\mathcal{P}, \mathcal{B})$ a $t$-$(v,k,\lambda)$ design, or simply {\em $t$-design}. The elements of $\mathcal{P}$ are
called
 {\em points}, and those of $\mathcal{B}$ are referred to as  {\em blocks}. We often define the number of blocks by $b$ and a $t$-design is simple when
 there is no
repeated blocks in $\mathcal{B}$. A $t$-design is called {\em symmetric} if $v=b$ and  trivial if $k=t$ or $k=v$. Hereafter we study only simple
$t$-design with
$t<k<v$.  We call $t$-$(v, k, \lambda)$ design a  {\em Steiner system} when $t\geq 2$ and $\lambda=1$. It is well known that a $t$-$(v,k,\lambda)$ design is also an $i$-$(v,k,\lambda_i)$ design for any $0\leq i< t,$ where
\begin{eqnarray}\label{condition}
\lambda_i=\frac{\lambda{{v-i}\choose {t-i}}}{ {{k-i}\choose {t-i}}}.
\end{eqnarray}

Combinatorial $t$-designs have very important applications in coding theory, cryptography, communications and statistics. There has been an interplay
between codes and $t$-designs for both directions. On  one hand, a linear code over any finite field can be derived from the incidence matrix of a
$t$-design and much progress has been
made and documented in \cite{Assmus-Key-1992-CUP, Ding-Ding-2015-IE, Tonchev-1998, Tonchev-2007}. On the other hand, linear and nonlinear codes might both
hold $t$-designs. Till  now,   $4$-designs and $5$-designs with fixed parameters were only derived from binary and ternary Golay codes. The largest $t$
for which the infinite families of $t$-designs could be obtained from linear codes is $t=3$.  In 2017, Ding and Li \cite{Ding-Li-2017-DM} obtained
infinite families of $2$-designs  from $p$-ary Hamming codes, ternary projective cyclic codes, binary codes with two zeros and their duals as well as the
infinite families of $3$-designs from the extended codes of these codes and RM codes. Afterwards, infinite families of $2$-designs and $3$-designs were
constructed
from a class of binary linear codes with five weights  \cite{Ding-2018-CC}. For other constructions of $t$-designs, for example, we refer to
\cite{Beth-Jungnickel-Lenz-1999-CUP, Colbowrn-Mafthon-2007-CB, MacWilliams-Sloane-1977-NHA, Reid-Rosa-2010-EJC}.

The rest of this paper is organized as follows. In Section \ref{section-2}, we introduce some notation, together with some preliminary results on
affine-invariant codes and $2$-designs, which will be used in subsequent sections. In Section \ref{theorem}, we present the weight distributions of a
class of linear codes derived from the dual of
the extended cyclic codes with two non-zeros. By proving that the derived codes are affine-invariant, we then obtain infinite families of $2$-designs and
explicitly determine their parameters. The proof of the main results are given in Section \ref{proofs}. Section \ref{section-5} concludes the paper.

\section{Preliminaries}\label{section-2}
In this section, we present some standard notation and basic facts on affine-invariant codes and $2$-designs.
\subsection{Some notation}
For convenience, we fix the following notation unless otherwise stated in this paper.
\begin{itemize}
\item $p$ is an odd prime and $p^*=(-1)^{\frac{p-1}{2}}p$.
\item $m\geq2$, $1\leq l\leq m-1$ are positive integers, $q=p^m$, $n=p^m-1$ and gcd$(m, l)=d>1$.
\item $\zeta_p=e^{2\pi \nu /p}$ is a primitive $p$-th root of unity, where $\nu=\sqrt{-1}$.
\item $\mathbb{F}_q$ denotes the finite field with $q$
elements and $\mathbb{F}^*_q=\mathbb{F}_q\backslash\{0\}$.
\item $\eta$ and $\eta'$ are the quadratic characters of $\mathbb{F}_q^*$ and $\mathbb{F}_p^*,$ respectively. We extend these quadratic characters by
    setting $\eta(0)=0$ and $\eta'(0)=0$.
\item $Tr$ denotes the trace function from $\mathbb{F}_q$ onto $\mathbb{F}_p.$
\end{itemize}

\subsection{The classical construction of $t$-designs from  affine-invariant codes}
Let $p$ be an odd prime, $m, l$ be positive integers and $(m, l)=d>1$. Let $\mathbb{F}_q$ denote the finite field with $q=p^m$
elements and $\mathbb{F}^*_q=\mathbb{F}_q\backslash\{0\}$. An $[n,k,d]$  {\em linear code} $\mathcal{C}$ over $\mathbb{F}_p$ is a $k$-dimensional subspace
of
$\mathbb{F}_p^n$ with minimum Hamming distance $d$, and is  {\em cyclic} provided that any cyclic shift of a codeword is another codeword of
$\mathcal{C}$. Any cyclic code $\mathcal{C}$ can be expressed as $\mathcal{C}= \langle g(x) \rangle,$ where $g(x)$ is monic and has the least degree. The
polynomial $g(x)$ is called the  {\em generator polynomial} and $h(x)=(x^n-1)/g(x)$ is referred to as the
{\em parity check polynomial} of $\mathcal{C}$. If the
generator polynomial $g(x)$ (resp. the parity  check polynomial $h(x)$) can be factored
into a product of $s$ irreducible polynomials over $\mathbb{F}_q$, then $\mathcal{C}$  is called a cyclic code with
 {\em $s$ zeros} (resp. {\em $s$ non-zeros}). The code with the generator polynomial $x^kh(x^{-1})/h(0)$ is called the {\em dual} of $\mathcal{C}$ and
 denoted by
$\mathcal{C}^{\bot}$, which is an $[n,n-k]$ code. The { \em extended code } of $\mathcal{C}$ is defined by
$$\overline\mathcal{C}=\{(c_0, c_1, \ldots, c_n) \in \mathbb{F}_p^{n+1}:(c_0, c_1, \ldots, c_{n-1}) \in \mathcal{C} ~\mathrm{with} ~\sum^n_{i=0}c_i=0\}.$$

For a codeword $\mathbf{c}\in \mathcal{C},$ the  {\em support} of $\mathbf{c}$ is defined by
$$Suppt(\mathbf{c})=\{0\leq i \leq n-1: c_i\neq 0\}.$$
Let $A_i$ be the number of codewords with Hamming weight $i$ in a code $\mathcal{C}$. The  {\em weight enumerator} of $\mathcal{C}$ is defined by
$$1+A_1z+A_2z^2+\ldots+A_nz^n,$$
and the sequence $(1, A_1, \ldots, A_n)$ is called the  {\em weight distribution} of  $\mathcal{C}.$ If $|\{1\leq i\leq n: A_i\neq 0\}|=w,$ then we say
$\mathcal{C}$ a  {\em $w$-weight} code.  Indexed the  coordinates
of a codeword in $\mathcal{C}$ by $(0, 1, 2, \ldots, n-1),$ then for each $i$ with $A_i\neq 0$, we denote  $\mathcal{B}_i$ by  the set of the supports of all codewords with weight
$i$. Let $\mathcal{P}=\{0, 1, \ldots, n-1\}.$ The pair $(\mathcal{P},\mathcal{B}_i)$ could give a
$t$-$(n,i,\lambda)$ design for some positive $\lambda$ \cite{Ding-Li-2017-DM} and the parameter $\lambda$ can be determined explicitly by Eq.(\ref{condition}). Such codes were used to construct $2$-designs and $3$-designs in
\cite{Assmus-Key-1992-CUP, Tonchev-1998, Tonchev-2007}.

We conclude this section by introducing some known results on affine-invariant codes
related to $2$-designs.

The set of coordinate permutations that map a code $\mathcal{C}$ to itself forms a group, which is called the  {\em permutation automorphism
group}  of
$\mathcal{C}$ and denoted by $PAut(\mathcal{C})$.
We define the {\em affine group} $GA_1(\mathbb{F}_q)$ by the set of all permutations
$$\sigma_{a,b}: x\mapsto  ax+b $$
of $\mathbb{F}_q$, where $a \in \mathbb{F}_q^*$ and  $b \in \mathbb{F}_q.$
An affine-invariant code is an extended cyclic code $\overline \mathcal{C}$ over $\mathbb{F}_p$ satisfying $GA_1(\mathbb{F}_q)\subseteq PAut(\overline
\mathcal{C})$ \cite{Huffman-Pless-2003-CUP}.

The $p$-adic expansion of each $s\in\mathcal{P}$ is given by
$$s=\sum^{m-1}_{i=0}s_ip^i,~ ~0\leq s_i\leq p-1 ,~0\leq i \leq m-1.$$
For any $r=\sum^{m-1}_{i=0}r_ip^i \in\mathcal{P}$,  we call that $r\preceq s$ if $r_i \leq s_i$ for all $0\leq i\leq m-1.$

For any integer $0\le j \le n-1$, the {\em $p$-cyclotomic coset} of $j$ modulo $n$ is defined by
$$C_j=\{jp^i \pmod n: 0 \le i \le \ell_j-1\}$$
where $\ell_j$ is the smallest positive integer such that $j\equiv jp^{\ell_j}\pmod n.$ The smallest integer in $C_j$ is called the {\em coset leader} of
$C_j$. Let $g(x)=\prod_j\prod_{i\in C_j}(x-\alpha^i)$, where $j$ runs through some coset leaders of the $p$-cyclotomic cosets $C_j$ modulo $n.$ The
 union of these $p$-cyclotomic cosets 
$T=\bigcup_jC_j$ is referred to as the {\em defining set} of $\mathcal{C}$.

For certain applications, it is important to know whether a
given extended primitive cyclic code $\overline\mathcal{C}$ is affine-invariant or not. The following Lemma \ref{Kasami-Lin-Peterson} presented by Kasami \emph{et al.} \cite{Kasami-Lin-Peterson-1968-IC} answered the questions by examining the defining set of the code. Lemma \ref{The dual of an affine-invariant code} proposed by Ding \cite{Ding-2019} in 2019 point out that one can get the new affine-invariant codes from the known ones.

\begin{lemma}\label{Kasami-Lin-Peterson}
 Let $\overline\mathcal{C}$ be an extended cyclic code of
length $p^m$ over $\mathbb{F}_p$ with defining set $\overline{T}$. The code $\overline\mathcal{C}$ is affine-invariant if and
only if whenever $s \in   \overline{T}$  then $r \in \overline{T}$ for all $r \in \mathcal{P}$  with $r \preceq s$.
\end{lemma}

\begin{lemma}\label{The dual of an affine-invariant code}
 The dual of an affine-invariant code $\overline\mathcal{C}$ over $\mathbb{F}_p$ of length $n+1$
is also affine-invariant.
\end{lemma}

Lemma \ref{The dual of an affine-invariant code} and the following Theorem \ref{2-design}  given by Ding in \cite{Ding-2019} are very powerful tools in
constructing $2$-designs from linear codes. We will employ them  to present our main results  in this paper.

\begin{theorem}\label{2-design}
 For each $i$ with $A_i \neq 0$ in an affine-invariant code $\overline \mathcal{C}$, the supports of the codewords of weight $i$ form a
$2$-design.
\end{theorem}

The following theorem presented by Ding in \cite{Ding-2019} reveals the relation of the codewords with the same support in a linear code  $\mathcal{C}$,
which will be used together with Eq.(\ref{condition}) to calculate the parameters of  $2$-designs later.

\begin{theorem}\label{design parameter}
 Let $\mathcal{C}$ be a linear code over $\mathbb{F}_p$ with minimum weight $d$. Let $w$ be the largest integer with $w\leq n$ satisfying
 $${w-\lfloor\frac{w+p-2}{p-1}\rfloor}<d.$$
Let $\mathbf{c_1}$ and $\mathbf{c_2}$ be two codewords of weight $i$ with $d\leq i\leq w$ and $Suppt(\mathbf{c_1})=Suppt(\mathbf{c_2}).$ Then
$\mathbf{c_1}=a\mathbf{c_2}$ for some $a\in \mathbb{F}_p^*$.
\end{theorem}

\section{The weight distribution of the codes and their $2$-designs}\label{theorem}
In this section, we present the main results, including the weight distribution and the underlying $2$-designs. The proofs of the main results will be
provided in the following section. We begin this section by the definition of the code.

Define
\begin{eqnarray}\label{code-2}
{\overline{{\mathcal{C}}^{\bot}}}^{\bot}:=\{ (Tr (ax^{p^l+1}+bx)_{x \in \mathbb{F}_q}+h ): a,b \in \mathbb{F}_q, h\in \mathbb{F}_p \},
\end{eqnarray}
where $\mathcal{C}$ is the  cyclic code of length $n$  and parity  check polynomial $h_1(x)h_2(x)$ where $h_1(x), h_2(x)$ are the minimal polynomials of
$\alpha, \alpha^{p^l+1}$ over $\mathbb{F}_p$, respectively, and $\alpha$ is a primitive element of $\mathbb{F}_q$. In particular, if $l=\frac{m}{2}$ we
have
$$
{\overline{{\mathcal{C}}^{\bot}}}^{\bot}:=\{ (Tr (ax^{p^{\frac{m}{2}}+1}+bx)_{x \in \mathbb{F}_q}+h ): a\in\mathbb{F}_{p^\frac{m}{2}}, b \in \mathbb{F}_q,
h\in \mathbb{F}_p
\}.
$$
Below, we consider only the case $gcd(m,l)>1$  since the case for $gcd(m,l)=1$ was discussed in \cite{Du-Wang-Tang}.

The following Theorems \ref{weight2}-\ref{parameter-3} are the main results of this paper.

\begin{theorem}\label{weight2}
Let $m\geq 2$ and $(m, l)=d>1$. The weight distribution of the code ${\overline{{\mathcal{C}}^{\bot}}}^{\bot}$ over $\mathbb{F}_p$ with length $n+1$ are
given in Table \ref{1} when both $m$ and $\frac{m}{d}$ are odd, in Table \ref{2} when $m$ is even and $\frac{m}{d}$ is odd, in Table \ref{3} when $\frac{m}{d}\geq 4$ is even
and in Table \ref{4} when $m=2d$, respectively. And the dimension is $dim({\overline{{\mathcal{C}}^{\bot}}}^{\bot})=\frac{3m}{2}+1$
when $m=2d$ and dim(${\overline{{\mathcal{C}}^{\bot}}}^{\bot})=2m+1$ otherwise.
\begin{table}
\begin{center}
\caption{The weight distribution of ${\overline{{\mathcal{C}}^{\bot}}}^{\bot}$ when both $m$ and $m/d$ are odd}\label{1}
\begin{tabular}{ll}
\hline\noalign{\smallskip}
Weight  &  Multiplicity   \\
\noalign{\smallskip}
\hline\noalign{\smallskip}
$0$  &  $1$ \\
$p^{m-1}(p-1)$  &  $ p(p^{m-1}+1)(p^m-1) $    \\
$ p^{m-1}(p-1)+p^{\frac{m-1}{2}}$  &  $ \frac{1}{2}p^m(p-1)(p^m-1)$     \\
$ p^{m-1}(p-1)-p^{\frac{m-1}{2}}$  &  $ \frac{1}{2}p^m(p-1)(p^m-1)$     \\
$p^m$ & $ p-1 $   \\
\noalign{\smallskip}
\hline
\end{tabular}
\end{center}
\end{table}

\begin{table}
\begin{center}
\caption{The weight distribution of ${\overline{{\mathcal{C}}^{\bot}}}^{\bot}$ when $m$ is even and $m/d$ is odd}\label{2}
\begin{tabular}{ll}
\hline\noalign{\smallskip}
Weight  &  Multiplicity   \\
\noalign{\smallskip}
\hline\noalign{\smallskip}
$0$  &  $1$ \\
$p^{m-1}(p-1)$  &  $ p(p^m-1) $    \\
$ p^{m-1}(p-1)+p^{\frac{m}{2}-1}$  &  $ \frac{1}{2}p^m(p-1)(p^m-1)$     \\
$ p^{m-1}(p-1)-p^{\frac{m}{2}-1}$  &  $ \frac{1}{2}p^m(p-1)(p^m-1)$     \\
$ p^{m-1}(p-1)+p^{\frac{m}{2}-1}(p-1)$  &  $ \frac{1}{2}p^m(p^m-1)$     \\
$ p^{m-1}(p-1)-p^{\frac{m}{2}-1}(p-1)$  &  $ \frac{1}{2}p^m(p^m-1)$     \\
$p^m$ & $ p-1 $  \\
\noalign{\smallskip}
\hline
\end{tabular}
\end{center}
\end{table}

\begin{table}
\begin{center}
\caption{The weight distribution of ${\overline{{\mathcal{C}}^{\bot}}}^{\bot}$ when $m/d\geq4$ is even}\label{3}
\begin{tabular}{ll}
\hline\noalign{\smallskip}
Weight  &  Multiplicity   \\
\noalign{\smallskip}
\hline\noalign{\smallskip}
$0$  &  $1$ \\
$p^{m-1}(p-1)$  &  $ p(p^{m-d}-p^{m-2d}+1)(p^m-1) $    \\
$ p^{m-1}(p-1)-(-1)^{\frac{m}{2d}}p^{\frac{m}{2}-1}(p-1)$  &  $ p^{m+d}(p^m-1)/(p^d+1)$     \\
$ p^{m-1}(p-1)+(-1)^{\frac{m}{2d}}p^{\frac{m}{2}-1}$  &  $ p^{m+d}(p-1)(p^m-1)/(p^d+1)$     \\
$ p^{m-1}(p-1)+(-1)^{\frac{m}{2d}}p^{\frac{m}{2}+d-1}(p-1)$  &  $ p^{m-2d}(p^m-1)/(p^d+1)$     \\
$ p^{m-1}(p-1)-(-1)^{\frac{m}{2d}}p^{\frac{m}{2}+d-1}$  &  $ p^{m-2d}(p-1)(p^m-1)/(p^d+1)$     \\
$p^m$ & $ p-1 $  \\
\noalign{\smallskip}
\hline
\end{tabular}
\end{center}
\end{table}

\begin{table}
\begin{center}
\caption{The weight distribution of ${\overline{{\mathcal{C}}^{\bot}}}^{\bot}$ when $m/d=2$}\label{4}
\begin{tabular}{ll}
\hline\noalign{\smallskip}
Weight  &  Multiplicity   \\
\noalign{\smallskip}
\hline\noalign{\smallskip}
$0$  &  $1$ \\
$p^{m-1}(p-1)$  &  $ p(p^m-1) $    \\
$ p^{m-1}(p-1)-p^{\frac{m}{2}-1}$  &  $ p^m(p^{\frac{m}{2}}-1)(p-1)$     \\
$ p^{m-1}(p-1)+p^{\frac{m}{2}-1}(p-1)$  &  $ p^m(p^{\frac{m}{2}}-1)$     \\
$p^m$ & $ p-1 $  \\
\noalign{\smallskip}
\hline
\end{tabular}
\end{center}
\end{table}
\end{theorem}

One can see that the code is four-weight if both $m$ and $\frac{m}{d}$ are odd or $m=2d$ and six-weight if $m$ is even and $\frac{m}{d}$ is odd or
$\frac{m}{d}\geq 4$ is even.

\begin{theorem}\label{$2-$design-1}
Let $m\geq 2, l$ be positive integers with $gcd(m,l)>1$. Then the supports of the codewords of weight $i$ with $A_i\neq0$ in
${\overline{{\mathcal{C}}^{\bot}}}^{\bot}$ form
a $2$-design.
\end{theorem}

\begin{theorem}\label{parameter-3}
Let $m\geq2, l$ be integers, $\mathcal{P}=\{0, 1, \ldots, p^m-1\}$ and $\mathcal{B}$ be the set of the supports of the codewords of
${\overline{{\mathcal{C}}^{\bot}}}^{\bot}$ with weight $i,$ where $A_i\neq 0.$ Then ${\overline{{\mathcal{C}}^{\bot}}}^{\bot}$ gives $2$-$(p^m, i,
\lambda)$ designs for the following pairs:

(1) if both $m$ and $m/d$ are odd,
\begin{itemize}
\item $(i, \lambda)=(p^m-p^{m-1}, (p^m-p^{m-1}-1)(p^{m-1}+1)).$
\item $(i, \lambda)=(p^m-p^{m-1}+p^{\frac{m-1}{2}},
    \frac{1}{2}p^{\frac{m-1}{2}}(p^{\frac{m+1}{2}}-p^{\frac{m-1}{2}}+1)(p^m-p^{m-1}+p^{\frac{m-1}{2}}-1)).$
\item $(i, \lambda)=(p^m-p^{m-1}-p^{\frac{m-1}{2}},
    \frac{1}{2}p^{\frac{m-1}{2}}(p^{\frac{m+1}{2}}-p^{\frac{m-1}{2}}-1)(p^m-p^{m-1}-p^{\frac{m-1}{2}}-1)).$
\end{itemize}

(2) if $m$ is even and $m/d$ is odd,
\begin{itemize}
\item $(i, \lambda)=(p^m-p^{m-1}, p^m-p^{m-1}-1).$
\item $(i, \lambda)=(p^{m-1}(p-1)+p^{\frac{m}{2}-1},
    \frac{1}{2}p^{\frac{m}{2}-1}(p^{\frac{m}{2}+1}-p^{\frac{m}{2}}+1)(p^m-p^{m-1}+p^{\frac{m}{2}-1}-1)).$
\item $(i, \lambda)=(p^{m-1}(p-1)-p^{\frac{m}{2}-1},
    \frac{1}{2}p^{\frac{m}{2}-1}(p^{\frac{m}{2}+1}-p^{\frac{m}{2}}-1)(p^m-p^{m-1}-p^{\frac{m}{2}-1}-1)).$
\item $(i, \lambda)=(p^{m-1}(p-1)+p^{\frac{m}{2}-1}(p-1),
    \frac{1}{2}p^{\frac{m}{2}-1}(p^{\frac{m}{2}}+1)(p^m-p^{m-1}+p^{\frac{m}{2}}-p^{\frac{m}{2}-1}-1)).$
\item $(i, \lambda)=(p^{m-1}(p-1)-p^{\frac{m}{2}-1}(p-1),
    \frac{1}{2}p^{\frac{m}{2}-1}(p^{\frac{m}{2}}-1)(p^m-p^{m-1}-p^{\frac{m}{2}}+p^{\frac{m}{2}-1}-1)).$
\end{itemize}

(3) if $m/d\geq4$ is even,
\begin{itemize}
\item $(i, \lambda)=(p^m-p^{m-1}, (p^m-p^{m-1}-1)(p^{m-d}-p^{m-2d}+1)).$
\item $(i, \lambda)=(p^{m-1}(p-1)+(-1)^{\frac{m}{2d}}p^{\frac{m}{2}-1},
    p^{\frac{m}{2}+d-1}(p^{\frac{m}{2}+1}-p^{\frac{m}{2}}+(-1)^{\frac{m}{2d}})(p^m-p^{m-1}+(-1)^{\frac{m}{2d}}p^{\frac{m}{2}-1}-1)/(p^d+1)).$
\item $(i, \lambda)=(p^{m-1}(p-1)-(-1)^{\frac{m}{2d}}p^{\frac{m}{2}-1}(p-1),
    p^{\frac{m}{2}+d-1}(p^{\frac{m}{2}}-(-1)^{\frac{m}{2d}})(p^m-p^{m-1}-(-1)^{\frac{m}{2d}}p^{\frac{m}{2}}+(-1)^{\frac{m}{2d}}p^{\frac{m}{2}-1}-1)
    /(p^d+1)).$
\item $(i, \lambda)=(p^{m-1}(p-1)+(-1)^{\frac{m}{2d}}p^{\frac{m}{2}+d-1}(p-1),
    p^{\frac{m}{2}-d-1}(p^{\frac{m}{2}-d}+(-1)^{\frac{m}{2d}})(p^m-p^{m-1}+(-1)^{\frac{m}{2d}}p^{\frac{m}{2}+d}-(-1)^{\frac{m}{2d}}p^{\frac{m}{2}+d
    -1}-1)/(p^d+1)).$
\item $(i, \lambda)=(p^{m-1}(p-1)-(-1)^{\frac{m}{2d}}p^{\frac{m}{2}+d-1},
    p^{\frac{m}{2}-d-1}(p^{\frac{m}{2}-d+1}-p^{\frac{m}{2}-d}-(-1)^{\frac{m}{2d}})(p^m-p^{m-1}-(-1)^{\frac{m}{2d}}p^{\frac{m}{2}+d-1}-1)/(p^d+1)).$
\end{itemize}

(4) if $m/d=2$,
\begin{itemize}
\item $(i, \lambda)=(p^m-p^{m-1}, p^m-p^{m-1}-1).$
\item $(i, \lambda)=(p^m-p^{m-1}-p^{\frac{m}{2}-1},
    p^{\frac{m}{2}-1}(p^{\frac{m}{2}+1}-p^{\frac{m}{2}}-1)(p^m-p^{m-1}-p^{\frac{m}{2}-1}-1)/(p^\frac{m}{2}+1)).$
\item $(i, \lambda)=(p^m-p^{m-1}+p^{\frac{m}{2}}-p^{\frac{m}{2}-1}, p^{\frac{m}{2}-1}(p^m-p^{m-1}+p^{\frac{m}{2}}-p^{\frac{m}{2}-1}-1)).$
\end{itemize}
\end{theorem}


The following examples from Magma program confirm the main results in Theorems \ref{weight2} and \ref{parameter-3}.

\begin{example}\label{example1}
If $(p, l, m)=(3, 2, 4)$, then the code ${\overline{{\mathcal{C}}^{\bot}}}^{\bot}$ has parameters $[81,7,51]$ and weight enumerator
$1+1296z^{51}+240z^{54}+648z^{60}+2z^{81},$ It gives $2$-$(81, i, \lambda)$ designs with the following pairs $(i, \lambda):$
$$(51, 255), (54, 53), (60, 177),$$
which confirmed the results in Theorems \ref{weight2} and \ref{parameter-3}.
\end{example}

\begin{example}\label{example2}
If $(p, l, m)=(3, 3, 6)$, then the code ${\overline{{\mathcal{C}}^{\bot}}}^{\bot}$ has parameters $[729,10,477]$ and weight enumerator
$1+37908z^{477}+2184z^{486}+18954z^{504}+2z^{729},$ which confirmed the results in Theorem \ref{weight2}.
\end{example}

\begin{example}\label{example3}
If $(p, l, m)=(3, 2, 6)$, then the code ${\overline{{\mathcal{C}}^{\bot}}}^{\bot}$ has parameters $[729,10,468]$ and weight enumerator
$1+265356z^{468}+530712z^{477}+2184z^{486}+530712z^{495}+265356z^{504}+2^{729},$ which confirmed the results in Theorem \ref{weight2}.
\end{example}

\section{Proofs of the main results}\label{proofs}
We begin this section by  stating some  notation and basic facts on exponential sums.
\subsection{Exponential sums}
An additive character of $\mathbb{F}_q$ is a non-zero  function $\chi$ from $\mathbb{F}_q$ to the set of complex numbers of absolute value $1$ such that
$\chi(x+y)=\chi(x)\chi(y)$ for any pair $(x,y) \in \mathbb{F}_q^2$. For each $u \in \mathbb{F}_q$, the function
$$\chi_u(v)=\zeta^{Tr(uv)},~v \in \mathbb{F}_q$$
 denotes an additive character of $\mathbb{F}_q$. $\chi_0$ is  the trivial additive character of
 $\mathbb{F}_q$ since $\chi_0(v)=1$ for all $v \in \mathbb{F}_q$. We call $\chi_1$ the canonical additive character of $\mathbb{F}_q$ and we have $\chi_u(x)=\chi_1(ux)$ for all $u\in\mathbb{F}_q$
 \cite{Lidl-Niederreeiter-1997-CUP}.

The Gauss sum $G(\eta', \chi'_1)$ over $\mathbb{F}_p$ is defined by $G(\eta', \chi'_1)=\sum\limits_{c \in \mathbb{F}_p^*}\eta'(v)\chi'_1(v)=\sum\limits_{v
\in \mathbb{F}_p}\eta'(v)\chi'_1(v)$, where $\chi'_1$ is the canonical additive characters of $\mathbb{F}_p$.

To determine the parameters of the code ${\overline{{\mathcal{C}}^{\bot}}}^{\bot}$ defined in Eq.(\ref{code-2}), we introduce the following function.
\begin{equation}\label{S(a.b)}
S(a,b)=\sum\limits_{x \in \mathbb{F}_q}\chi(ax^{p^l+1}+bx), \quad a,b \in \mathbb{F}_q,
\end{equation}
if $m\ne 2l$, and otherwise
$$S(a,b)=\sum\limits_{x \in \mathbb{F}_q}\chi(ax^{p^{\frac{m}{2}}+1}+bx), \quad a\in \mathbb{F}_{p^{\frac{m}{2}}}, ~b \in \mathbb{F}_q.$$

The following Lemmas \ref{Gauss}-\ref{solution} are essential to determine the values of Eq.(\ref{S(a.b)}).

\begin{lemma}\label{Gauss}
\cite{Lidl-Niederreeiter-1997-CUP} With the notation above, we have
$$G(\eta', \chi'_1)={\sqrt{(-1)}^{(\frac{p-1}{2})^2}}\sqrt{p}=\sqrt{p^*}.$$
\end{lemma}

\begin{lemma}\label{Qu-Fea}
\cite{Ding-Ding-2015-IE} For each $y \in \mathbb{F}_p^*$,  $\eta(y)=1$ if $m\geq2$ is even, and $\eta(y)=\eta'(y)$ if  $m\geq1$ is odd.
\end{lemma}

\begin{lemma}\label{polynomial}
\cite{Coulter-1998-AA} Let $f(x)=a^{p^l}x^{p^{2l}}+ax \in \mathbb{F}_q[x]$, gcd$(m, l)=d$ and $b \in \mathbb{F}_q.$ There are three cases.

$(1)$ If $m/d$ is odd, then $f(x)$ is a permutation polynomial over $\mathbb{F}_q$ and
$$S(a,b)=\sqrt{p^*}^m\eta(a)\chi_1(-ax_{a,b}^{p^l+1}).$$

$(2)$ If $m/d$ is even and $a^{\frac {q-1}{p^d+1}}\neq (-1)^{\frac{m}{2d}},$ then $f(x)$ is a permutation polynomial over $\mathbb{F}_q$ and
$$S(a,b)=(-1)^{\frac{m}{2d}}p^{m/2}\chi_1(-ax^{p^l+1}_{a,b}).$$

$(3)$ If $m/d$ is even and $a^{\frac {q-1}{p^d+1}}=(-1)^{\frac{m}{2d}}$, then $f(x)$ is not a permutation polynomial over $\mathbb{F}_q$. We have
$S(a,b)=0$ when the equation $f(x)=-b^{p^l}$ is unsolvable, and
$$S(a,b)=-(-1)^{\frac{m}{2d}}p^{m/2+d}\chi_1(-ax^{p^l+1}_{a,b})$$
otherwise.
Particularly,
\begin{eqnarray*}
S(a,0)=\left\{
\begin{array}{ll}
\sqrt{p^*}^m\eta(a)   & \mathrm{if}\,\ m/d ~is ~odd,\\
(-1)^{\frac{m}{2d}}p^{\frac{m}{2}}    & \mathrm{if}\,\ m/d ~is ~ even~and~ a^{\frac{q-1}{p^d+1}} \neq (-1)^{\frac{m}{2d}}, \\
(-1)^{\frac{m}{2d}+1}p^{\frac{m}{2}+d}    & \mathrm{if}\,\ m/d ~ is ~even~and~  a^{\frac{q-1}{p^d+1}} = (-1)^{\frac{m}{2d}}.
\end{array}
\right.
\end{eqnarray*}
\end{lemma}
Notice that $x_{a,b}^{p^l+1}$ is a solution to the equation $f(x)=-b^{p^l}$. Moreover, $x_{a,b}^{p^l+1}$ is the unique solution when $f(x)$ is a
permutation polynomial over $\mathbb{F}_q$.

\begin{lemma}\label{solution}
\cite{Coulter-1991-M} For $m$  even and gcd$(m, l)=d$, the equation $a^{p^l}x^{p^{2l}}+ax=0$ is solvable for $x \in \mathbb{F}_q^*$ if and only if
$$a^{\frac{q-1}{p^d+1}}=(-1)^{\frac{m}{2d}}.$$
In such cases there are $p^{2d}-1$ non-zero solutions.
\end{lemma}

%

\subsection{Proofs of the main results}

Now we are ready to give the proofs of our main results. First we present the proof of the weight distribution of
${\overline{{\mathcal{C}}^{\bot}}}^{\bot}$
which is given in Eq.(\ref{code-2}).

$\mathbf{Proof~ of ~Theorem~ \ref{weight2}.}$
For each non-zero codeword $\mathbf{c}(a,b,h)=(c_0,\ldots,c_n)$ in ${\overline{{\mathcal{C}}^{\bot}}}^{\bot},$ the Hamming weight of $\mathbf{c}(a,b,h)$
is
\begin{equation}\label{weight formula-2}
w_H(\mathbf{c}(a,b,h))=p^m-T(a,b,h),
\end{equation}
where $$T(a,b,h)=|\{x:Tr(ax^{p^l+1}+bx)+h=0,~ x,a,b \in \mathbb{F}_q, h \in \mathbb{F}_p\}|$$
if $l\neq \frac{m}{2}$, and otherwise
$$T(a,b,h)=|\{x:Tr(ax^{p^{\frac{m}{2}}+1}+bx)+h=0,~ a\in \mathbb{F}_{p^{\frac{m}{2}}},  x,b \in \mathbb{F}_q, h \in \mathbb{F}_p\}|.$$

Then
\begin{eqnarray*}
T(a,b,h)&=&{\frac{1}{p}}\sum\limits_{y \in \mathbb{F}_p}\sum\limits_{x \in \mathbb{F}_q}\zeta_p^{yTr(ax^{p^l+1}+bx)+hy}\\
&=&p^{m-1}+{\frac{1}{p}}\sum\limits_{y\in \mathbb{F}_p^*}\zeta_p^{hy}\sum\limits_{x \in \mathbb{F}_q}\zeta_p^{yTr(ax^{p^l+1}+bx)}.
\end{eqnarray*}
If $a=b=h=0$, then $\mathbf{c}(a,b,h)$ is the zero codeword.\\
If $a=b=0, h\neq 0$, then  $T(a,b,h)=p^{m-1}+p^{m-1}\sum\limits_{y \in \mathbb{F}_p^*}\zeta_p^{hy}=0.$\\
If $a=0, b\neq 0$, then $ T(a,b,h)=p^{m-1}+{\frac{1}{p}}\sum\limits_{y \in \mathbb{F}_p^*}\zeta_p^{hy}\sum\limits_{x\in
\mathbb{F}_q}\zeta_p^{yTr(bx)}=p^{m-1}.$\\
If $a\neq 0$, then from the definition of $S(a,b)$, we have
\begin{eqnarray*}
T(a,b,h)
&=&p^{m-1}+{\frac{1}{p}}\sum\limits_{y \in \mathbb{F}_p^*}\zeta_p^{hy}S(ay,by).
\end{eqnarray*}
According to Lemma \ref{polynomial}, the calculation of $T(a,b,h)$ can be divide into the following cases.

(1) When both $m$ and $m/d$ are odd, by Lemma \ref{polynomial}, we have
\begin{eqnarray*}
T(a,b,h)=\left\{
\begin{array}{ll}
p^{m-1} & \mathrm{if}\,\ h=Tr(ax_{ay,by}^{p^l+1}),\\
p^{m-1}+p^{\frac{m-1}{2}}(-1)^\frac{(p-1)(m+1)}{4}\eta(a)\eta'(h-Tr(ax_{a,b}^{p^l+1})) & \mathrm{if}\,\ h\neq Tr(ax_{ay,by}^{p^l+1}).\\
\end{array}
\right.
\end{eqnarray*}
Moreover, we find that  $ T(a,b,h)=p^{m-1}$ appears $p^m(p^m-1)$ times, and both $$T(a,b,h)=p^{m-1}+p^{\frac{m-1}{2}}
~\mathrm{and}~ T(a,b,h)=p^{m-1}-p^{\frac{m-1}{2}}$$  appear $$\frac{p^m(p-1)(p^m-1)}{2}$$ times, respectively, from the properties of the quadratic characters $\eta$ of $\mathbb{F}_q^*$.

(2) When $m$ is even and $m/d$ is odd, by Lemma \ref{polynomial}, we get
\begin{eqnarray*}
T(a,b,h)=\left\{
\begin{array}{ll}
p^{m-1}-p^{\frac{m}{2}-1}(-1)^\frac{(p-1)m}{4}\eta(a) & \mathrm{if}\,\ h\neq Tr(ax_{ay,by}^{p^l+1}),\\
p^{m-1}+p^{\frac{m}{2}-1}(-1)^\frac{(p-1)m}{4}\eta(a)(p-1) & \mathrm{if}\,\ h= Tr(ax_{ay,by}^{p^l+1}).\\
\end{array}
\right.
\end{eqnarray*}
More precisely, by the properties of the quadratic characters $\eta$, we have  that both $$T(a,b,h)=p^{m-1}+p^{\frac{m}{2}-1}~ \mathrm{and} ~ T(a,b,h)=p^{m-1}-p^{\frac{m}{2}-1}$$  appear
$$\frac{p^m(p-1)(p^m-1)}{2}$$ times, and both $$T(a,b,h)=p^{m-1}+p^{\frac{m}{2}-1}(p-1)~ \mathrm{and}~ T(a,b,h)=p^{m-1}-p^{\frac{m}{2}-1}(p-1)$$  appear
$$\frac{p^m(p^m-1)}{2}$$ times, respectively.

(3) When $m/d$ is even and $a^{\frac{q-1}{p^d+1}}=(-1)^{\frac{m}{2d}},$ from Lemma \ref{polynomial}, we get
\begin{eqnarray*}
T(a,b,h)=\left\{
\begin{array}{ll}
p^{m-1}-(-1)^{\frac{m}{2d}}p^{\frac{m}{2}+d-1}(p-1) & \mathrm{if}\,\ f(x)=-b^{p^l}\,  \mathrm{is~ solvable}\, \mathrm{and}\, \\
&h=Tr(ax_{ay,by}^{p^l+1}),\, \\
p^{m-1}+(-1)^{\frac{m}{2d}}p^{\frac{m}{2}+d-1} & \mathrm{if}\,\ f(x)=-b^{p^l}\,  \mathrm{is ~solvable},\, b\neq 0\, \mathrm{and}\\
&\ h\neq
Tr(ax_{ay,by}^{p^l+1})\, \mathrm{or}\, b=0\,\mathrm{and}\, h\neq 0,\\
p^{m-1}  & \mathrm{if}\, f(x)=-b^{p^l}\,  \mathrm{is ~ no~ solvable}.
\end{array}
\right.
\end{eqnarray*}
It follows from  Lemma \ref{solution} that  there are $\frac{q-1}{p^d+1}$ elements $a\in \mathbb{F}_q^*$ satisfying $a^{\frac{q-1}{p^d+1}}=(-1)^{\frac{m}{2d}},$ and
$p^{m-2d}$ elements $b\in \mathbb{F}_q$ such that $f(x)=-b^{p^l}$ is solvable. Therefore,
$$T(a,b,h)=p^{m-1}+(-1)^{\frac{m}{2d}+1}p^{\frac{m}{2}+d-1}(p-1)$$
appears ${\frac{p^{m-2d}(p^m-1)}{p^d+1}}$ times,
$$T(a,b,h)=p^{m-1}+(-1)^{\frac{m}{2d}}p^{\frac{m}{2}+d-1}$$
appears ${\frac{p^{m-2d}(p^m-1)(p-1)}{p^d+1}}$ times, and $$T(a,b,h)=p^{m-1}$$
appears $(p^m-1)p^{m-2d+1}(p^d-1)$
times.

(4)  When $m/d \ge 4$ is even and $a^{\frac{q-1}{p^d+1}}\neq (-1)^{\frac{m}{2d}},$ by Lemma \ref{polynomial}, we have
\begin{eqnarray*}
T(a,b,h)=\left\{
\begin{array}{ll}
p^{m-1}+(-1)^{\frac{m}{2d}}p^{\frac{m}{2}-1}(p-1) & \mathrm{if}\,\ h=Tr(ax_{ay,by}^{p^l+1}),\\
p^{m-1}-(-1)^{\frac{m}{2d}}p^{\frac{m}{2}-1} & \mathrm{if}\,\  h\neq Tr(ax_{ay,by}^{p^l+1})\,.
\end{array}
\right.
\end{eqnarray*}
Particularly, there exist $p^m-1-{\frac{p^m-1}{p^d+1}}={\frac{p^d(p^m-1)}{p^d+1}}$
elements $a\in \mathbb{F}_q^*$ such that
$$a^{\frac{q-1}{p^d+1}}\neq (-1)^{\frac{m}{2d}}.$$
Then
$$T(a,b,h)=p^{m-1}+(-1)^{\frac{m}{2d}}p^{\frac{m}{2}-1}(p-1)$$ appears
${\frac{p^{m+d}(p^m-1)}{p^d+1}}$ times, and
 $$T(a,b,h)=p^{m-1}-(-1)^{\frac{m}{2d}}p^{\frac{m}{2}-1}$$
 appears ${\frac{p^{m+d}(p-1)(p^m-1)}{p^d+1}}$
times.

(5)  When $m=2d$, it is clear that  $a^{p^{\frac{m}{2}}-1}\neq -1$ for all $a\in \mathbb{F}_{p^{\frac{m}{2}}}^*$,  Similarly, by Lemma \ref{polynomial}, we have
\begin{eqnarray*}
T(a,b,h)&=&|\{x:Tr(ax^{p^{\frac{m}{2}}+1}+bx)+h=0,~ a\in \mathbb{F}_{p^{\frac{m}{2}}},  x,b \in \mathbb{F}_q, h \in \mathbb{F}_p\}|\\
&=&\left\{
\begin{array}{ll}
p^{m-1}-p^{\frac{m}{2}-1}(p-1) & \mathrm{if}\,\ h=Tr(ax_{ay,by}^{p^l+1}),\\
p^{m-1}+p^{\frac{m}{2}-1} & \mathrm{if}\,\  h\neq Tr(ax_{ay,by}^{p^l+1})\,.
\end{array}
\right.
\end{eqnarray*}
More precisely, we get
$$T(a,b,h)=p^{m-1}-p^{\frac{m}{2}-1}(p-1)$$ appears
$(p^{\frac{m}{2}}-1)p^m $ times, and
 $$T(a,b,h)=p^{m-1}+p^{\frac{m}{2}-1}$$
 appears $(p^{\frac{m}{2}}-1)(p-1)p^m $
times.

From all the discussions above, we complete the proof of the theorem.\hfill$\square$

Now we prove that $\overline{{\mathcal{C}}^{\bot}}$ is affine invariant.
\begin{lemma}\label{affine-invariant}
The extended code $\overline{{\mathcal{C}}^{\bot}}$ is affine-invariant.
\end{lemma}

\emph{Proof.} We will prove the conclusion with Lemma \ref{Kasami-Lin-Peterson}.
The defining set $T$ of the cyclic code ${\mathcal{C}}^{\bot}$ is $T =C_1\cup C_{p^l+1}$. Since $0 \not \in T$, the defining set $\overline{T}$ of
$\overline{{\mathcal{C}}^{\bot}}$ is given by $\overline{T} = C_1\cup C_{p^l+1} \cup \{0\}$.
Let $s \in \overline{T} $ and $r \in \mathcal{P}$. Assume that $r \preceq s$. We need to prove that $r \in \overline{T}$ by Lemma
\ref{Kasami-Lin-Peterson}.

If $r=0,$ then obviously $r\in \overline{T}$. Consider now the case $r>0$. If $s \in  C_1$, then the Hamming weight $wt(s) = 1.$ As $r \preceq s$, $wt(r)
= 1.$ Consequently, $r \in C_1 \subset  \overline{T}.$  If $s  \in C_{p^l+1}$, then the Hamming weight $wt(s) = 2.$ As $r \preceq s$, either $wt(r) = 1$
or $r = s.$ In both cases, $r \in  \overline{T}.$ The desired conclusion then follows
from Lemma \ref{Kasami-Lin-Peterson}.
Thus we complete the proof of the lemma.\hfill $\square$

Now we are ready to prove Theorems \ref{$2-$design-1} and \ref{parameter-3}.

\textbf{Proof of Theorem \ref{$2-$design-1}.}

By Lemmas \ref{affine-invariant} and \ref{The dual of an affine-invariant code}, we find that ${\overline{{\mathcal{C}}^{\bot}}}^{\bot}$ is
affine-invariant.  Thus  the supports of the codewords of weight $i$ with $A_i\neq0$ in ${\overline{{\mathcal{C}}^{\bot}}}^{\bot}$ give  $2$-designs from
Theorem \ref{2-design}. \hfill $\square$

\textbf{Proof of Theorem \ref{parameter-3}.}

By Theorem \ref{design parameter}, one can prove that the number of supports of all codewords with weight $i\neq 0$ in the code
${\overline{{\mathcal{C}}^{\bot}}}^{\bot}$ is equal to $A_i/(p-1)$ for each $i,$ where $A_i$ is given in Tables \ref{1}-\ref{4}. Then the desired
conclusions follow from Theorem \ref{$2-$design-1} and  Eq.(\ref{condition}). Thus, we complete the proof. \hfill$\square$

\section{Concluding remarks}\label{section-5}
In this paper, first we determined the weight distribution of a class of linear codes derived from the dual of extended cyclic codes with two non-zeros.
Using the properties of affine-invariant codes, we then found that the code ${\overline{{\mathcal{C}}^{\bot}}}^{\bot}$ holds $2$-designs. By computing
certain exponential sums, we explicitly determined their
parameters. The total number of the $2$-designs we obtained is exponential in $m$.

\end{document}